\documentclass{amsart}

\usepackage[dvips]{graphicx}
\usepackage{amssymb}
\usepackage{amsfonts}
\usepackage{verbatim}
\usepackage{pstricks, pst-tree, pst-node}
\usepackage{diagrams}
\usepackage{mathrsfs}
\usepackage[all]{xy}

\newtheorem{thm}{Theorem}
\newtheorem{defn}[thm]{Definition}
\newtheorem{rmk}[thm]{Remark}

\newcommand{\ot}{\otimes}

\title{Rectifying Partial Algebras Over Operads of Complexes}
\author{Scott O. Wilson}
\date{October 4, 2010}

\begin{document}


\begin{abstract}
In \cite{KM} Kriz and May introduced partial algebras over an operad.
In this paper we prove that, in the category of chain complexes, partial
algebras can be functorially replaced by quasi-isomorphic algebras.
In particular, partial algebras contain all of the important homological and
homotopical information that genuine algebras do. Applying this result 
to McClure's partial algebra in \cite{JM} shows that the 
chains of a PL-manifold are quasi-isomorphic to an $E_\infty$-algebra.
\end{abstract}

\maketitle

\section{Introduction}

In \cite{KM} Kriz and May introduced partial algebras over an operad and
proved that, in the category of simplicial modules, such partial
algebras are quasi-isomorphic to genuine algebras.
It was left as an open question whether or not such a result also holds in
the category of chain complexes. This is important given the recent work
of McClure, showing that the chains of a PL-manifold
form a partial algebra \cite{JM}.

In this paper we prove that, in the category of chain complexes, partial
algebras can be functorially replaced by quasi-isomorphic algebras.
In particular, partial algebras contain all of the important homological and
homotopical information that genuine algebras do. Applying this result 
to McClure's partial algebra in \cite{JM} shows that the 
chains of a PL-manifold are quasi-isomorphic to an $\mathcal{E}_\infty$-algebra.
We describe further applications and consequences as well.

I thank Jim McClure, Dennis Sullivan and the referee for their comments and 
suggestions.

\section{Preliminaries}
In this section we give some definitions.
We work in the category of flat chain complexes over some Dedekind ring $R$.
By a simplicial complex we mean a simplicial object in the category of 
complexes (so such an object is bi-graded).

\begin{defn} An operad (of complexes over $R$) 
is a collection of complexes $\mathcal{O}(j)$ over $R$, $j \geq 0$, 
together with a unit map $\eta: R \to \mathcal{O}(1)$, 
an action of the symmetric
group $\Sigma_j$ on $\mathcal{O}(j)$ for each $j$, and chain maps
$$
\gamma : \mathcal{O}(k) \otimes \mathcal{O}(j_1) \otimes \cdots \otimes 
\mathcal{O}(j_k) \to \mathcal{O}(j_1 + \cdots + j_k)
$$
for all $k \geq 1$ and $j_i \geq 0$. The maps $\gamma$ are required to be
associative, equivariant with respect to the $\Sigma$-actions, and unital
with respect to the unit $\eta$. See \cite{KM}.
\end{defn}

Intuitively, the component $\mathcal{O}(j)$ encodes operations with $j$ inputs and one output. 
The maps $\gamma$ determine the composition of operations. 
Morphisms of operads are defined naturally.

An algebra over an operad $\mathcal{O}$ is a complex $X$ with chain
maps 
$$
\mathcal{O}(j) \otimes_{R[\Sigma_j]} X^{\otimes_j} \to X
$$ 
sending the operad unit to the identity map of $X$, and satisfying the obvious axiom codifying
an action with respect to operad composition. Here $\Sigma_j$
acts on  $\mathcal{O}(j) \otimes_{R[\Sigma_j]} A^{\otimes_j}$ by $\sigma$
on the left and $\sigma^{-1}$ on the right.

To define a partial algebra
we first introduce the notion of a \emph{domain} on which an operad may
partially act. This first appeared in \cite{KM}.

\begin{defn} \label{defn:domain}
A domain in a complex $X$ is a collection of subcomplexes
\[
i_j: X_j \to X^{\otimes_j}
\]
satisfying the following:
\begin{enumerate}
\item $X_1 = X$.
\item For all $j = j_1 + \dots + j_k$, $X_j$ is a $\Sigma_j$-invariant 
subcomplex of $X_{j_1} \ot \dots \ot X_{j_k}$, making the following diagram 
commute
\[
\xymatrix{
 X_j \ar[r]^-{i_{j_1, \ldots ,j_k}} \ar[d]_-{i_j} & X_{j_1} \ot \dots \ot X_{j_k} \ar[d]^{i_{j_1} \ot \dots \ot i_{j_k}}  \\
 X^{\ot_j} \ar[r]^-{\cong} & X^{\ot_{j_1}} \ot \dots \ot X^{\ot_{j_k}}
 }
\]
\item The inclusion map $i_j: X_j \hookrightarrow X^{\otimes_j}$ is a  quasi-isomorphism.
\end{enumerate}
\end{defn}

We remark that our flatness assumption and condition $(3)$ imply that the inclusions $i_{j_1, \ldots ,j_k}$
in condition $(2)$ are quasi-isomorphisms.

A morphism $f$ of domains $\{X_j\}$ and  $\{Y_j\}$ is a collection of chain maps $f_j : X_j \to Y_j$
such that each map $f_j$ equals the restriction of ${f_1}^{\ot_j}$ to $X_j$.
We say $f$ is a quasi-isomorphism if each $f_j$ is a quasi-isomorphism. 
It follows from our flatness 
assumption, and the diagram
\[
\xymatrix{
 X_j \ar[r]^-{f_j} \ar[d] & Y_j \ar[d] \\
 X^{\ot_j} \ar[r]^-{{f_1}^{\ot_j}} & Y^{\ot_j}
 }
\]
that if $f_1$ is a quasi-isomorphism, then each $f_j$ is a quasi-isomorphism.

\begin{rmk} \label{rmk:RL}
There is a functor $L$ from domains to complexes taking 
$\{X_j\}$ to $X_1 = X$. There is also a functor $R$ from complexes to 
domains taking $X$ to the domain $X_j = X^{\otimes_j}$, and $LR=id$.
\end{rmk}

\begin{defn} \label{defn:partialalg}
Let $\mathcal{O}$ be an operad. A partial algebra over the operad 
$\mathcal{O}$ is a domain $\{X_j\}$ in a complex $X$ and a collection of 
chain maps
$$
\Theta_j: \mathcal{O}(j) \otimes_{R[\Sigma_j]} X_j \to X
$$
satisfying the following:
\begin{enumerate}
\item The operad unit acts as the identity: $\Theta_1 \circ (\eta \ot id_X) = id_X$.
\item For all $j = j_1 + \dots + j_k$, the maps 
$$
\Theta_{j_1 , \ldots , j_k }: \mathcal{O}(j_1) \otimes \dots \otimes \mathcal{O}(j_k) \otimes_{R[\Sigma_j]} 
X_j \to X^{\ot_k}
$$ 
given by including $X_j$ into $X_{j_1} \otimes \dots \otimes X_{j_k}$, applying the shuffle, 
and then applying $\Theta_{j_1} \otimes \dots \otimes \Theta_{j_k} $, must factor through $X_k$.
\item The maps $\Theta_j$ describe an action with respect to the operad 
composition. Namely, for all $j = j_1 + \dots + j_k$,
$$
\Theta_j \circ (\gamma \otimes id_{X_j}) = \Theta_k \circ ( id_{\mathcal{O}(k)} \otimes \Theta_{j_1 , \ldots , j_k } )
$$
as maps from $\mathcal{O}(k) \otimes \mathcal{O}(j_1) \otimes \cdots \otimes 
\mathcal{O}(j_k) \otimes X_j$ to $X$.
\end{enumerate}
\end{defn}

A morphism of partial algebras over an operad is a morphism of domains that
commutes with the partial actions. We say a morphism of partial algebras is a 
quasi-isomorphism if it is a quasi-isomorphism of domains.

\begin{rmk} \label{paex}
An algebra over an operad is a partial algebra where the
domain $\{X_j\}$ is given by $X_j = RX = X^{\otimes_j}$.
\end{rmk}

\begin{figure}
\begin{center}
\begin{pspicture}(0,0)(0,2)
\pstree[treemode=U,levelsep=25pt]{\Tp}
{
  \pstree{\Tp}
  { 
    \Tp 
    \Tp 
    \Tp
    \Tp
  }
}
\end{pspicture}
\end{center}
\caption{A tree as an operation.} \label{tree}
\end{figure}
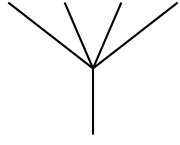

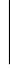
\begin{figure}
\begin{center}
\begin{pspicture}(0,0)(0,1)
\pstree[treemode=U,levelsep=25pt]{\Tp}{\Tp}
\end{pspicture}
\end{center}
\caption{The unit element.} \label{unit}
\end{figure}

We now give a diagrammatic description of operads and their algebras.
We represent elements of $\mathcal{O}(k)$ by trees with $k$ inputs, as in 
Figure~\ref{tree}, and the unit in $\mathcal{O}(1)$ as in Figure \ref{unit}.
Implicit in this are the various structures of an operad over complexes:
addition, the differential and the $\Sigma_j$-actions.

\begin{figure}
\begin{center}
\begin{pspicture}(0,0)(4,4)
\rput(2,2){
\pstree[treemode=U,levelsep=25pt]{\Tp}
{
  \pstree{\Tp}
  {
    \pstree{\Tcircle[linecolor=white]{}}
    {
      \pstree{\Tp}
      { \Tp 
	\Tp 
	\Tp
      } 
      
    }
    \Tp
    \Tp
    \pstree{\Tcircle[linecolor=white]{}}
    {
       \pstree{\Tp}
       {
	 \Tp 
	 \Tp 
	 \Tp
       }
    }
  }
}
}
\rput(2,3){$\cdots$}
\end{pspicture}
\end{center}
\caption{A generator of $\mathcal{O}(k) \otimes_{R[\Sigma_k]} \mathcal{O}(j_1)
\cdots \mathcal{O}(j_k)$.} \label{trees}
\end{figure}
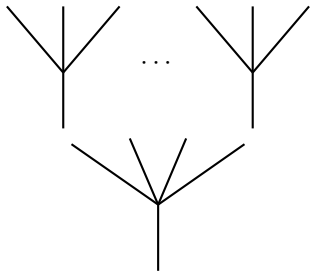

We represent a generator of $\mathcal{O}(k) \otimes_{R[ \Sigma_k ]} \mathcal{O}(j_1) \otimes 
\cdots \otimes \mathcal{O}(j_k)$ by a collection of $k$ trees, as in 
Figure \ref{trees}, where we have left spaces between trees to indicate this 
is a tensor product of elements of $\mathcal{O}$. Again various
structures are implicit. In particular, the symmetric group acts on the bottom tree, and also
by permuting the tensor factor of trees on top.

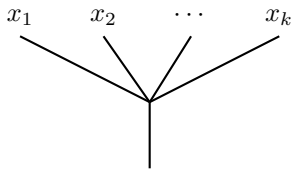
\begin{figure}
\begin{center}
\begin{pspicture}(0,0)(0,2.5)
\pstree[treemode=U,levelsep=25pt]{\Tp}
{
  \pstree[labelsep=2pt]{\Tp}
  {
    \Tp~{$x_1$}
    \Tp~{$x_2$}
    \Tp~{$\cdots$}
    \Tp~{$x_k$}
  }
}
\end{pspicture}
\end{center}
\caption{A generator of $\mathcal{O}(k) \otimes_{R[\Sigma_k]} X^{\otimes_k}$.} \label{OboxX}
\end{figure}

We represent a generator of 
\[
\mathcal{O} \boxtimes X = \sum_{k \geq 0} \mathcal{O}(k) \otimes_{R[\Sigma_k]} X^{\otimes_k}
\]
by a diagram consisting
of a tree labeled by elements of $X$, as in Figure \ref{OboxX},
 where $x_1 \otimes \cdots \otimes x_k \in X^{\otimes_k}$. 
There are still implicit
notions of addition and differential, as well as the symmetric group
actions. In particular, this picture is equivalent to the one obtained
by acting on the tree by $\sigma$ and on (the tensor product of)
the $k$ labeling elements by $\sigma^{-1}$, for all $\sigma \in \Sigma_k$.

There is a categorical interpretation of operads as monads which allows one to make use of the two sided bar construction.
The constructions and proofs in the next section involve minor variations of this bar construction.
Rather than appealing to this categorical construction abstractly, we will
unravel it explicitly in the case of partial algebras. We do this
because it makes our work more transparent, and secondly, because
it may be used to give a picture for the bar construction in more
general situations. 

Finally, some terminology. By the total complex of a simplicial complex $X_{q,k}$, with simplicial grading $q$ and
complex grading $k$, we mean the complex whose degree $n$ is $\sum_{p+k=n} X_{q,p}$ and whose differential is equal to the sum of the simplicial differential  $\sum (-1)^i \partial_i$ and $(-1)^q$ times the complex differential. 
Similarly, maps of simplicial complexes can be added along total degrees to give maps of total complexes.

\section{Main Result} \label{sec:algresult}

In this section we prove the following:

\begin{thm} \label{thm:convert} 
Let $A$ be a flat complex and $\mathcal{O} = \bigoplus_{k \geq 0} 
\mathcal{O}(k)$ be an operad of 
complexes such that each $\mathcal{O}(k)$ is a projective 
$R[\Sigma_k]$-module. There is a
functor $W$ that assigns to any partial $\mathcal{O}$-algebra $A_*$ an 
$\mathcal{O}$-algebra $WA_*$ such that $A_*$ and $WA_*$ are 
quasi-isomorphic as partial $\mathcal{O}$-algebras.
\end{thm}

Let us first give an outline for the proof. We will construct a diagram 
\begin{equation*} \label{diagram}
\begin{diagram}
A_* & \pile{\lTo^{\varphi} \\ \rTo_{\eta}} & B_* 
& \rTo^{\delta} & W_*
\end{diagram}
\end{equation*}
where $B_*$ is a partial $\mathcal{O}$-algebra and $W_*$ is an $\mathcal{O}$-algebra (and therefore also 
a partial $\mathcal{O}$-algebra, by Remark \ref{paex}). The maps $\eta$ and $\varphi$ are quasi-isomorphisms of domains. Moreover,
$\varphi$ is a morphism of partial $\mathcal{O}$-algebras, and therefore a quasi-isomorphism of 
partial $\mathcal{O}$-algebras. Finally, $\delta$ is a quasi-isomorphism of partial $\mathcal{O}$-algebras.
The constructions of $B_*$ and $W_*$ from $A_*$ will be seen to be natural, and the assignment $A_* \mapsto W_*$
will be the desired functor in the statement of the theorem.
 
The rest of this section is divided into subsections which complete the steps in this outline. Several of
the techniques used appear in \cite{KM}. 

\subsection{Definition of the complex $B$}

First we define a simplicial complex associated to the partial algebra $A$, whose $q$-simplicies are denoted by $B_q$. 
This first appeared in \cite{KM} (Definition 3.2, Example 4.2). The reader may note that this simplicial complex is a minor variation on the two sided bar construction $B( \mathcal{O}, \mathcal{O},C)$ where $C$ is an $ \mathcal{O}$-algebra.  

Let $A$ be a partial $\mathcal{O}$-algebra with domain $A_*$ having inclusions $i_j : A_j \to A^{\ot_j}$ and
$i_{\alpha_1, \ldots , \alpha_k} : A_\alpha \to A_{\alpha_1} \ot \dots \ot A_{\alpha_k}$, where $\alpha = \alpha_1 + \dots + \alpha_k$.
We let $B_0$ be the following subcomplex of $\mathcal{O} \boxtimes A$, induced by the domain $A_*$ and 
the inclusion maps $id \ot i_j$:
\[
B_0 = \bigoplus_{j \geq 0} \mathcal{O}(j) \otimes_{R[\Sigma_j]} A_j
\] 
Next we define 
\begin{align*}
B_1 &= \bigoplus_{\stackrel{k \geq 0}{j_1, \ldots, j_k \geq 0}} \mathcal{O}(k) \ot_{R[\Sigma_k]} \mathcal{O}(j_1) 
\ot \dots \ot \mathcal{O}(j_k) \otimes_{R[\Sigma_j]} A_{j}
\end{align*}
where $j = j_1 + \dots + j_k$.
In other words, $B_1$ is the natural the subcomplex of $\mathcal{O} \boxtimes \mathcal{O} \boxtimes A$
induced by the domain $A_*$.

Now, we consider general $q$. The complex
\[
\underbrace{\mathcal{O} \boxtimes \dots \boxtimes \mathcal{O}}_{q+1} \boxtimes A
\]
is naturally given as a direct sum of tensor products 
\[
\mathcal{O}(n_1) \ot \dots \ot \mathcal{O}(n_m) \ot A_{\alpha_1} \ot \dots \ot A_{\alpha_k}.
\]
By our flatness assumption, each such summand has a subcomplex 
$\mathcal{O}(n_1) \ot \dots \ot \mathcal{O}(n_m) \ot A_\alpha$, given by tensoring the given inclusions 
$i_{\alpha_1, \ldots ,\alpha_k} : A_\alpha \to A_{\alpha_1} \ot \dots \ot A_{\alpha_k}$ with the identity map on copies of $\mathcal{O}$. We let $B_q$ be the direct sum of these subcomplexes.

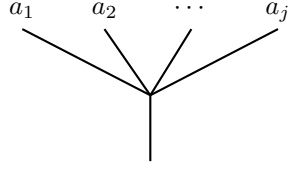
\begin{figure}
\begin{center}
\begin{pspicture}(0,0)(0,2.5)
\pstree[treemode=U,levelsep=25pt]{\Tp}
{
  \pstree[labelsep=2pt]{\Tp}
  {
    \Tp~{$a_1$}
    \Tp~{$a_2$}
    \Tp~{$\cdots$}
    \Tp~{$a_j$}
  }
}
\end{pspicture}
\end{center}
\caption{A generator of $B_0$ where $a_1 \ot \dots \ot a_j \in A_j$.} \label{labeledtree2}
\end{figure}

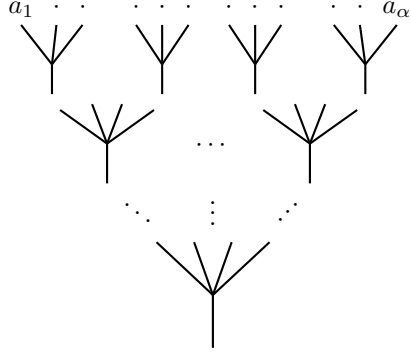
\begin{figure}
\begin{center}
\begin{pspicture}(0,0)(6,5)
\rput(3,.8){
\pstree[treemode=U,levelsep=20pt,treesep=.5]{\Tp}
{
  \pstree[labelsep=2pt]{\Tp}
  {
    \Tp
    \Tp
    \Tp
    \Tp
  }
}
}

\rput(2,2){$\ddots$}
\rput(3,2){$\vdots$}
\rput(4,1.9){$\cdot$}
\rput(3.9,1.8){$\cdot$}
\rput(4.1,2){$\cdot$}

\rput(1.5,3.5){
\pstree[treemode=U,levelsep=15pt,treesep=.4]{\Tp}
{
  \pstree{\Tp}
  {
    \pstree{\Tcircle[linecolor=white]{}}
    {
      \pstree[labelsep=1pt,treesep=.25]{\Tp}
      { \Tp~{$a_1$} 
	\Tp~{$\cdot$} 
	\Tp~{$\cdot$}
      } 
      
    }
    \Tp
    \Tp
    \pstree{\Tcircle[linecolor=white]{}}
    {
       \pstree[labelsep=1pt,treesep=.25]{\Tp}
       {
	 \Tp~{$\cdot$} 
	 \Tp~{$\cdot$} 
	 \Tp~{$\cdot$}
       }
    }
  }
}
}

\rput(3,2.8){$\cdots$}

\rput(4.4,3.5){
\pstree[treemode=U,levelsep=15pt,treesep=.4]{\Tp}
{
  \pstree{\Tp}
  {
    \pstree{\Tcircle[linecolor=white]{}}
    {
      \pstree[labelsep=1pt,treesep=.25]{\Tp}
      { \Tp~{$\cdot$}
	\Tp~{$\cdot$}
	\Tp~{$\cdot$}
      } 
      
    }
    \Tp
    \Tp
    \pstree{\Tcircle[linecolor=white]{}}
    {
       \pstree[labelsep=1pt,treesep=.25]{\Tp}
       {
	 \Tp~{$\cdot$} 
	 \Tp~{$\cdot$} 
	 \Tp~{$a_{\alpha}$}
       }
    }
  }
}
}

\end{pspicture}
\end{center}
\caption{A generator of $B_q$, with $a_1 \ot \dots \ot a_{\alpha} \in A_\alpha$.} \label{labeledtree3}
\end{figure}

There is a simple diagrammatic description of these complexes $B_q$. For example a generator of $B_0$ can be represented as in Figure \ref{labeledtree2}.  Also, a generator of $B_q$ can be represented, as in Figure \ref{labeledtree3}, as a stacking of trees of height $q+1$, with elements $a_1 \ot \dots \ot a_{\alpha} \in A_\alpha$.

The complexes $B_q$ form a simplicial complex with face and degeneracy maps defined on $B_q$ in exactly the same way as they are for the two sided bar construction $B(\mathcal{O}, \mathcal{O}, C)$ where $C$ is an $\mathcal{O}$-algebra. 
Namely, the zeroth face operator $\partial_0$ is given by the partial action of $\mathcal{O}$ on $A$, and for $0<i\leq q$, the $i^{th}$
face operator is given by the operad composition $\mathcal{O} \boxtimes \mathcal{O} \to \mathcal{O}$ in the $i^{th}$
$\boxtimes$-factor. It follows from the definition of $B_q$, and condition $(2)$ of the definition of partial algebras, that $\partial_i :B_q \to B_{q-1}$ for all $0 \leq i \leq q$.

Similarly, the $i^{th}$ face operator $s_i: B_q \to B_{q+1}$ is induced by the operadic unit $R \to \mathcal{O}$
in the $i^{th}$ $\boxtimes$-factor. The proof that this forms a simplicial set is exactly the same as the proof for the 
bar construction, and appears in \cite{KM} using the language of monads.

It is easy to visualize $B_q$ and its simplicial structure in terms of our diagrams.
As in Figure~\ref{labeledtree3}, 
let us refer to top row of trees at the $1^{st}$, the next below the $2^{nd}$,
etc. The $0^{th}$ face operator of this simplicial object is given by 
evaluating the elements of $A$ on the $1^{st}$ row of trees using the 
partial algebra structure of $A_*$ over $\mathcal{O}$.
For $1 \leq i \leq q$, the $i^{th}$ face operator is given by composing the
$i^{th}$ and $(i+1)^{st}$ rows of this diagram using the operad structure.
The $0^{th}$ degeneracy operator of this simplicial object is given by 
inserting a row of units of  $\mathcal{O}$ between the elements of $A$ and 
the first row of trees. For $1 \leq i \leq q$,
the $i^{th}$ degeneracy operator is given by inserting a row of units of 
$\mathcal{O}$ between the $i^{th}$ and $(i+1)^{st}$ rows of this diagram.

This completes our definition and description of the simplicial complex $B$. We will use the same notation, $B$, to denote total complex associated to the simplicial complex $B$.

\subsection{Domain in $B$} \label{Bdomain}

We now define a domain in the total complex $B$, and denote the $j^{th}$ subcomplex of $B^{\ot_j}$
by $B_j$, or the entire domain simply by $B_*$. For $j=1$ we let $B_1 = B$. 

For $j>1$ we note that, again as before, $B^{\ot_j}$ can be written as a sum over terms which are given by a tensor product of copies of the operad components $\mathcal{O}(n_i)$ and of the subcomplexes $A_{\alpha_i}$.  Again by the definition of the domain $A_*$ and our flatness assumption, each summand has a subspace induced by the given inclusions 
$i_{\alpha_1, \ldots , \alpha_k} : A_\alpha \to  A_{\alpha_1} \ot \dots \ot A_{\alpha_k} $, where $\alpha = \alpha_1 + \dots + \alpha_k$. 
We let $B_j$ be the sum of these subspaces of $B^{\ot_j}$, and we denote the induced inclusion by 
$I_j :B_j \to B^{\ot_j}$.
It follows from conditions $(1)$ and $(2)$ in the definition of partial algebra that $B_j$ is a subcomplex of $B^{\ot_j}$.

\begin{figure}
\begin{center}
\begin{pspicture}(0,0)(6,3)
\rput(0,1.5){
\pstree[treemode=U,levelsep=25pt]{\Tp}
{
  \pstree[labelsep=2pt]{\Tp}
  {
    \Tp~{$a_1$}
    \Tp~{$a_2$}
    \Tp~{$\cdot$}
    \Tp~{$\cdot$}
  }
}
}
\rput(3,1.5){$\cdots$}
\rput(6,1.5){
\pstree[treemode=U,levelsep=25pt]{\Tp}
{
  \pstree[labelsep=2pt]{\Tp}
  {
    \Tp~{$\cdot$}
    \Tp~{$\cdot$}
    \Tp~{$a_{\alpha-1}$}
    \Tp~{$a_\alpha$}
  }
}
}
\end{pspicture}
\end{center}
\caption{A generator in the domain $B_j$, with $a_1 \ot \dots \ot a_\alpha \in A_\alpha$.} \label{B0j}
\end{figure}

There is a diagrammatic description of $B_j$ given as follows. For simplicial degree $q=0$ and arbitrary complex degree, we can view a generator as $j$ trees labeled by $a_1 \ot \dots \ot a_\alpha \in A_\alpha$, as in Figure \ref{B0j}.

More generally, we can picture a generator of 
$B_j$ that is contained in $B_q^{\ot_j}$, for some fixed $q$,
as in Figure \ref{Bqk}, by a stacking of trees each of height $q+1$ labeled on top by elements of $A$
such that $a_1 \otimes \cdots \otimes a_\alpha \in A_\alpha$. 
In the most general case, we can picture a generator as a collection of $j$ stackings of trees, all of various heights, labeled on top by elements of $A$ such that $a_1 \otimes \cdots \otimes a_\alpha \in A_\alpha$ (not shown).

\begin{figure}
\begin{center}
\begin{pspicture}(0,0)(12,5)
\rput(3,.8){
\pstree[treemode=U,levelsep=20pt,treesep=.5]{\Tp}
{
  \pstree[labelsep=2pt]{\Tp}
  {
    \Tp
    \Tp
    \Tp
    \Tp
  }
}
}

\rput(2,2){$\ddots$}
\rput(3,2){$\vdots$}
\rput(4,1.9){$\cdot$}
\rput(3.9,1.8){$\cdot$}
\rput(4.1,2){$\cdot$}

\rput(1.5,3.5){
\pstree[treemode=U,levelsep=15pt,treesep=.4]{\Tp}
{
  \pstree{\Tp}
  {
    \pstree{\Tcircle[linecolor=white]{}}
    {
      \pstree[labelsep=1pt,treesep=.25]{\Tp}
      { \Tp~{$a_1$} 
	\Tp~{$\cdot$} 
	\Tp~{$\cdot$}
      } 
      
    }
    \Tp
    \Tp
    \pstree{\Tcircle[linecolor=white]{}}
    {
       \pstree[labelsep=1pt,treesep=.25]{\Tp}
       {
	 \Tp~{$\cdot$} 
	 \Tp~{$\cdot$} 
	 \Tp~{$\cdot$}
       }
    }
  }
}
}

\rput(3,2.8){$\cdots$}

\rput(4.4,3.5){
\pstree[treemode=U,levelsep=15pt,treesep=.4]{\Tp}
{
  \pstree{\Tp}
  {
    \pstree{\Tcircle[linecolor=white]{}}
    {
      \pstree[labelsep=1pt,treesep=.25]{\Tp}
      { \Tp~{$\cdot$}
	\Tp~{$\cdot$}
	\Tp~{$\cdot$}
      } 
      
    }
    \Tp
    \Tp
    \pstree{\Tcircle[linecolor=white]{}}
    {
       \pstree[labelsep=1pt,treesep=.25]{\Tp}
       {
	 \Tp~{$\cdot$} 
	 \Tp~{$\cdot$} 
	 \Tp~{$\cdot$}
       }
    }
  }
}
}

\rput(1.6,3.9){$\cdots$}
\rput(4.4,3.9){$\cdots$}

\rput(6.5,2){$\cdots$}

\rput(10,.8){
\pstree[treemode=U,levelsep=20pt,treesep=.5]{\Tp}
{
  \pstree[labelsep=2pt]{\Tp}
  {
    \Tp
    \Tp
    \Tp
    \Tp
  }
}
}

\rput(9,2){$\ddots$}
\rput(10,2){$\vdots$}
\rput(11,1.9){$\cdot$}
\rput(10.9,1.8){$\cdot$}
\rput(11.1,2){$\cdot$}

\rput(8.5,3.5){
\pstree[treemode=U,levelsep=15pt,treesep=.4]{\Tp}
{
  \pstree{\Tp}
  {
    \pstree{\Tcircle[linecolor=white]{}}
    {
      \pstree[labelsep=1pt,treesep=.25]{\Tp}
      { \Tp~{$\cdot$} 
	\Tp~{$\cdot$} 
	\Tp~{$\cdot$}
      } 
      
    }
    \Tp
    \Tp
    \pstree{\Tcircle[linecolor=white]{}}
    {
       \pstree[labelsep=1pt,treesep=.25]{\Tp}
       {
	 \Tp~{$\cdot$} 
	 \Tp~{$\cdot$} 
	 \Tp~{$\cdot$}
       }
    }
  }
}
}

\rput(10,2.8){$\cdots$}

\rput(11.4,3.5){
\pstree[treemode=U,levelsep=15pt,treesep=.4]{\Tp}
{
  \pstree{\Tp}
  {
    \pstree{\Tcircle[linecolor=white]{}}
    {
      \pstree[labelsep=1pt,treesep=.25]{\Tp}
      { \Tp~{$\cdot$}
	\Tp~{$\cdot$}
	\Tp~{$\cdot$}
      } 
      
    }
    \Tp
    \Tp
    \pstree{\Tcircle[linecolor=white]{}}
    {
       \pstree[labelsep=1pt,treesep=.25]{\Tp}
       {
	 \Tp~{$\cdot$} 
	 \Tp~{$\cdot$} 
	 \Tp~{$a_\alpha$}
       }
    }
  }
}
}

\rput(8.6,3.9){$\cdots$}
\rput(11.3,3.9){$\cdots$}
\end{pspicture}
\end{center}
\caption{A generator of $B_j$. Here with $a_1 \ot \dots \ot a_\alpha \in A_\alpha$.} \label{Bqk}
\end{figure} 

It remains to show that $B_j$ is a domain in $B$. First, $B_1 = B$ by definition, and the $\Sigma_j$-equivariance 
of $I_j: B_j \to B^{\ot_j}$ follows from the $\Sigma$-equivariance in the domain $A_*$ and the definition of $B_j$. Secondly, 
the inclusion map $I_j : B_j \to B^{\ot_j}$ satisfies the factoring condition $(2)$ in 
Definition \ref{defn:domain} since it is induced by the inclusions $i_*$ of the domain $A_*$ which satisfy this condition.  Lastly, we claim the inclusions $I_j : B_j \to B^{\ot_j}$ are quasi-isomorphism since they are
induced by the inclusions in the domain $A_*$, that are quasi-isomorphisms. 

To prove this, note for each $j > 1$ there is a spectral sequence for each of the bi-complexes $B_j$ and $B^{\ot_j}$ and a morphism between them given by $I_j$. 
The induced map on the first page is an isomorphism since these pages are the homology with respect
to the differentials on $\mathcal{O}$ and $A_*$, the inclusions $i_{\alpha_1, \ldots , \alpha_k}$ are all quasi-isomorphisms, and tensoring with the projective $R [ \Sigma_{n_i}]$-module $\mathcal{O}(n_i)$ preserves quasi-isomorphisms. These spectral sequences converge to the homology of their total complexes, since they are bounded, and therefore the induced map on the homology is an isomorphism since it was an isomorphism on the first page.

\subsection{Quasi-isomorphism of domains $B_*$ and $A_*$} \label{BAqi}

We first construct quasi-isomorphisms $\eta: A \to B$ and $\varphi: B \to A$. This was first done in \cite{KM} (Example 4.2), and is word for word the same as the proof that the usual bar construction is a resolution, so we will be brief. 

Let $\underline{A} $ denote the constant simplicial object with $A$ in each simplicial degree and all face and degeneracy maps given by the identity. There are canonical maps $\gamma: A \to \underline{A}$ and $\epsilon: \underline{A} \to A$
which are quasi-isomorphisms.

Next, we construct a chain equivalence of $\underline{A}$ and $B$. There is an inclusion $\psi:\underline{A}_q \to B_q$ 
of simplicial complexes given by 
\[
\psi(a) = \underbrace{u \ot \dots \ot u}_{q+1} \ot a 
\]
where $u \in \mathcal{O}(1)$ is the operad unit. Next, we have a map $\tau: B \to \underline{A} $ of simplicial complexes
given in each simplicial degree by evaluation using the partial action of $\mathcal{O}$ on $A$. By conditions $(2)$ 
and $(3)$ of the definition of partial algebra, this is well defined. It is immediate to check that $\tau \circ \psi = id$. 
Also, there is an explicit simplicial homotopy $h$ such that $\psi \circ \tau$ is homotopic to the identity. It is induced by the simplicial face operators defined in the previous section, as in the usual proof that the bar construction is a resolution, see \cite{KM}.  Therefore, we have quasi-isomorphisms on the total complexes
\begin{equation*} 
\begin{diagram}
A & \pile{\lTo^{\epsilon} \\ \rTo_{\gamma}} & \underline{A} & \pile{\lTo^{\tau} \\ \rTo_{\psi}} & B
\end{diagram}
\end{equation*}
and it follows that $\varphi = \epsilon \circ \tau: B \to A$ and $\eta = \psi \circ \gamma$ are quasi-isomorphisms. We will refer to the map $\varphi$ as the evaluation map since, for $x \in B$ of simplicial degree zero,  $\varphi(x) \in A$ is given by the partial action of $\mathcal{O}$ on $A$ (while for higher simplicial degrees the map is zero).

Now our goal is to show $\varphi$ and $\eta$ each induce a quasi-isomorphism of the domains $B_*$ and $A_*$.
Note that by condition $(2)$ of Definition \ref{defn:partialalg} and our definition of $B_j$, the restriction of the evaluation map
$\varphi^{\ot_j} : B^{\ot_j} \to A^{\ot_j}$ to $B_j$ factors through $A_j$. Similarly,
the restriction of the inclusion map $\eta^{\ot_j}:  A^{\ot_j} \to B^{\ot_j}$ to $A_j$ factors through $B_j$.
Thus we have a diagram
\[
\xymatrix{
B_j \ar[r]^-{I_j}  \ar@/_1pc/[d]_-{\varphi_j} & B^{\ot_j} \ar@/^1pc/[d]^-{\varphi^{\ot_j} } \\
A_j  \ar[r]_{i_j} \ar[u]_-{\eta_j}& A^{\ot_j} \ar[u]^-{\eta^{\ot_j}}
}
\]
where the square commutes starting from $A_j$ or $B_j$. By the remark after definition \ref{defn:domain}, 
$\varphi_j$ and $\eta_j$ are quasi-isomorphisms. Namely the top, bottom and right vertical maps are quasi-isomorphisms, so $\varphi_j: B_j \to A_j$ is also a quasi-isomorphism. For shorthand, we denote this quasi-isomorphism of domains $\varphi_j : B_j \to A_j$
by $\varphi$, as in the outline above. Similarly, $\eta_j$ is a quasi-isomorphism and gives a quasi-isomorphism $\eta$
of domains.
After defining the partial $\mathcal{O}$-algebra on $B_*$ below, we will see that $\varphi$ is a map of partial 
$\mathcal{O}$-algebras.

\subsection{$B_*$ as a partial $\mathcal{O}$-algebra} \label{paB}

In \cite{KM} the authors define, from a partial $\mathcal{O}$-algebra $A_*$, a simplicial partial algebra. 
We will recall it's definition here, as it will be used to define a partial $\mathcal{O}$-algebra on the domain $B_*$.

For each $q,j$ there as a natural subcomplex
\[
C_{q,j} \to \underbrace{B_q \ot \dots \ot B_q}_j
\]
given as before by the subspaces induced by the domain $A_*$ in each summand of the right hand side. In particular,
$C_{q,1} = B_q$. Moreover, there are chain maps
\[
\mathcal{O}(j) \ot_{R[ \Sigma_j]} C_{q,j} \to C_{q,1}
\]
given by the left action of $\mathcal{O}(j)$ on itself, and it is immediate to check this defines a simplicial partial 
$\mathcal{O}$-algebra. For each $j$ we let $C_j$ be induced total complex and note that $C_1 = B$.

Algebraic structures defined in terms of simplicial maps can be used to define algebraic 
structure on the chain level, by using the shuffle map (\cite{M}, Appendix, or \cite{KM}). 
Recall, that for simplicial complexes $X_{q,r}, Y_{p,s}$ the shuffle map 
\[
g:T(X)_{q+r} \otimes T(Y)_{p+s} \to T(X \otimes Y)_{q+p+r+s} 
\]
is defined on these total complexes by
\[
g(a \otimes b) = \sum_{(u,v)} \pm(s_{\nu_q} \cdots s_{\nu_1} a \otimes
s_{\mu_p} \cdots s_{\mu_1} b)
\]
where $s_*$ are the degeneracy operators, the sum is over all shuffles $\nu_1 < \cdots < \nu_q$ and 
$\mu_1 < \cdots < \mu_p$ of $\{0,1,\cdots,p+q+1\}$, and the sign is determined by the signature of the corresponding permutation of  $\{0,1,\cdots,p+q+1\}$. It is important to note that $g$ is commutative, associative and unital, and we denote the iterates of the shuffle map also by $g$.

In our case of interest, we have for each $q = q_1 + \dots + q_k$, the shuffle map
\[
g: B_{q_1} \ot \dots \ot B_{q_j} \to B_q^{\ot_j}
\]
We let $G_j$, with domain $B^{\ot_j}$, denote the sum of these shuffle maps. 
It is immediate from the definitions of $B_j$ and $C_j$ that the restriction of $G_j$ to $B_j$ 
factors through $C_j$, so that there is a well define map $\Sigma_j$-equivariant map $B_j \to B^{\ot_j} \to C_j$ 
given by the inclusion $I_j$ followed by the shuffle map. Thus we can define $\Theta_j$ to be the composition
\[
\mathcal{O}(j) \ot_{R[ \Sigma_j]} B_j \to \mathcal{O}(j) \ot_{R[ \Sigma_j]} C_j \to C_1 = B
\]

In words, this map is given by applying the shuffle product (to obtain an element of correct total degree) followed by the left action of $\mathcal{O}$ on itself. This indeed defines a partial $\mathcal{O}$-algebra on the domain $B_*$.
Property $(2)$ of Definition \ref{defn:partialalg} follows from (and in fact motived) the definition of $B_j$, while  property $(3)$ of Definition \ref{defn:partialalg} follows from the properties of $g$ mentioned above.

Diagrammatically, this action corresponds to inserting operadic units, according to the shuffle map, into a collection of $j$ stackings of trees, followed by composing the trees at the bottom of the diagram with a generator of $\mathcal{O}$.

We close this subsection by noting that the evaluation map $\varphi$ from subsection \ref{BAqi} is a map of partial $\mathcal{O}$-algebras, 
by properties $(1)$ and $(2)$ of partial $\mathcal{O}$-algebras. By the previous subsection, it follows that $\varphi$ is a quasi-isomorphism of partial $\mathcal{O}$-algebras.

\subsection{The $\mathcal{O}$-algebra $W_*$} \label{W}

In this section we define an $\mathcal{O}$-algebra $W_*$.
First, we define a simplicial complex $W$.
We let $W_0 = \mathcal{O} \boxtimes A$ and for $q \geq 1$ we let 
\[
W_q = \bigoplus_{k \geq 0} \mathcal{O}(k) \otimes_{R[\Sigma_k]} B_{q-1}^{\ot_k}
\]

The face and degeneracy operators of $W$ are defined in the same way as for $B$, using the 
the operad partial action for the zeroth face, the operad composition for the other face operations,
 and operad unit for the degeneracies. We will use the same notation $W$ for the induced total complex.

We note that there is a canonical map $\delta: B \to W$ again induced by the quasi-isomorphisms $i_j$ and
$i_{\alpha_1, \ldots , \alpha_k}$ of the domain
$A_*$ and the identity on all $\mathcal{O}$ tensor factors.

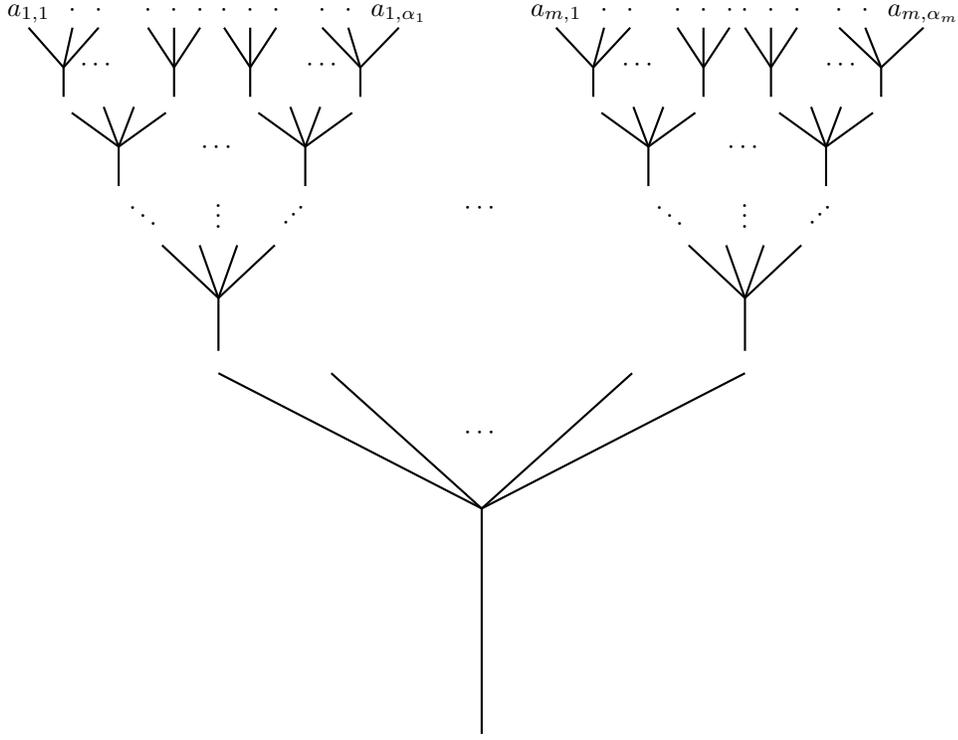
\begin{figure}
\begin{center}
\begin{pspicture}(0,-5)(12,5)
\rput(3,.8){
\pstree[treemode=U,levelsep=20pt,treesep=.5]{\Tp}
{
  \pstree[labelsep=2pt]{\Tp}
  {
    \Tp
    \Tp
    \Tp
    \Tp
  }
}
}

\rput(2,2){$\ddots$}
\rput(3,2){$\vdots$}
\rput(4,1.9){$\cdot$}
\rput(3.9,1.8){$\cdot$}
\rput(4.1,2){$\cdot$}

\rput(1.5,3.5){
\pstree[treemode=U,levelsep=15pt,treesep=.4]{\Tp}
{
  \pstree{\Tp}
  {
    \pstree{\Tcircle[linecolor=white]{}}
    {
      \pstree[labelsep=1pt,treesep=.25]{\Tp}
      { \Tp~{$a_{1,1}$} 
	\Tp~{$\cdot$} 
	\Tp~{$\cdot$}
      } 
      
    }
    \Tp
    \Tp
    \pstree{\Tcircle[linecolor=white]{}}
    {
       \pstree[labelsep=1pt,treesep=.25]{\Tp}
       {
	 \Tp~{$\cdot$} 
	 \Tp~{$\cdot$} 
	 \Tp~{$\cdot$}
       }
    }
  }
}
}

\rput(3,2.8){$\cdots$}

\rput(4.4,3.5){
\pstree[treemode=U,levelsep=15pt,treesep=.4]{\Tp}
{
  \pstree{\Tp}
  {
    \pstree{\Tcircle[linecolor=white]{}}
    {
      \pstree[labelsep=1pt,treesep=.25]{\Tp}
      { \Tp~{$\cdot$}
	\Tp~{$\cdot$}
	\Tp~{$\cdot$}
      } 
      
    }
    \Tp
    \Tp
    \pstree{\Tcircle[linecolor=white]{}}
    {
       \pstree[labelsep=1pt,treesep=.25]{\Tp}
       {
	 \Tp~{$\cdot$} 
	 \Tp~{$\cdot$} 
	 \Tp~{$a_{1,\alpha_1}$}
       }
    }
  }
}
}

\rput(1.4,3.9){$\cdots$}
\rput(4.4,3.9){$\cdots$}

\rput(6.5,2){$\cdots$}

\rput(10,.8){
\pstree[treemode=U,levelsep=20pt,treesep=.5]{\Tp}
{
  \pstree[labelsep=2pt]{\Tp}
  {
    \Tp
    \Tp
    \Tp
    \Tp
  }
}
}

\rput(9,2){$\ddots$}
\rput(10,2){$\vdots$}
\rput(11,1.9){$\cdot$}
\rput(10.9,1.8){$\cdot$}
\rput(11.1,2){$\cdot$}

\rput(8.5,3.5){
\pstree[treemode=U,levelsep=15pt,treesep=.4]{\Tp}
{
  \pstree{\Tp}
  {
    \pstree{\Tcircle[linecolor=white]{}}
    {
      \pstree[labelsep=1pt,treesep=.25]{\Tp}
      { \Tp~{$a_{m,1}$} 
	\Tp~{$\cdot$} 
	\Tp~{$\cdot$}
      } 
      
    }
    \Tp
    \Tp
    \pstree{\Tcircle[linecolor=white]{}}
    {
       \pstree[labelsep=1pt,treesep=.25]{\Tp}
       {
	 \Tp~{$\cdot$} 
	 \Tp~{$\cdot$} 
	 \Tp~{$\cdot$}
       }
    }
  }
}
}

\rput(10,2.8){$\cdots$}

\rput(11.4,3.5){
\pstree[treemode=U,levelsep=15pt,treesep=.4]{\Tp}
{
  \pstree{\Tp}
  {
    \pstree{\Tcircle[linecolor=white]{}}
    {
      \pstree[labelsep=1pt,treesep=.25]{\Tp}
      { \Tp~{$\cdot$}
	\Tp~{$\cdot$}
	\Tp~{$\cdot$}
      } 
      
    }
    \Tp
    \Tp
    \pstree{\Tcircle[linecolor=white]{}}
    {
       \pstree[labelsep=1pt,treesep=.25]{\Tp}
       {
	 \Tp~{$\cdot$} 
	 \Tp~{$\cdot$} 
	 \Tp~{$a_{m,\alpha_m}$}
       }
    }
  }
}
}

\rput(8.6,3.9){$\cdots$}
\rput(11.3,3.9){$\cdots$}

\psline(6.5,-5)(6.5,-2)
\psline(6.5,-2)(3,-.2)
\psline(6.5,-2)(4.5,-.2)
\rput(6.5,-1){$\cdots$}
\psline(6.5,-2)(8.5,-.2)
\psline(6.5,-2)(10,-.2)

\end{pspicture}
\end{center}
\caption{A generator of $W_q$. Here $a_{i,1} \otimes \cdots \otimes a_{i,\alpha_i} \in A_{\alpha_i}$, for all $1 \leq i \leq m$.} \label{fig2}
\end{figure}

There is a simple diagrammatic description of $W_q$. As in Figure \ref{fig2}, a generator of $W_q$ may be 
represented by a stacking of trees $q+1$ high, labeled on top by elements of $A$ 
satisfying the following property: for each tree in the $q^{th}$ (second to bottom) row, the elements 
$a_{i,1}, \cdots ,a_{i, \alpha_i}$ of $A$ ``lying above'' this tree
satisfy $a_{i,1} \otimes \cdots \otimes a_{i,\alpha_i} \in A_{\alpha_i}$.
Diagrammatically, the simplicial structure can be view as it is for $B$.

By pre-composition with the shuffle map, followed by the left action of $\mathcal{O}$ on itself, we can define a map
\[
\mathcal{O}(j) \ot_{R[ \Sigma_j]} W^{\ot_j} \to W
\] 
similar to subsection \ref{paB}. 
The new point to check is that the image of this composition does in fact land in $W$, but this follows from (and in fact motivates) the definition of $W$. The other properties follow as in the previous subsection  \ref{paB}, so that $W$ is an $\mathcal{O}$-algebra. Said another way, by Remark \ref{paex}, $W_*$ is a partial $\mathcal{O}$-algebra with domain $W_j = W^{\otimes_j}$ for all $j \geq 1$.

\subsection{A quasi-isomorphism of partial $\mathcal{O}$-algebras $B_*$ and $W_*$.}

As noted in the previous subsection \ref{W}, that there is a canonical map of complexes $\delta: B \to W$ induced by the quasi-isomorphisms 
of the domain $A$. It follows from the definition of $B_j$, and property $(2)$ of the domain $A_*$, that 
the restriction of $\delta^{\ot_j}$ to $B_j$ factors through $W^{\ot_j}$, giving a map of domains $\delta_j : B_j \to W_j = W^{\ot_j}$. 

This map is a quasi-isomorphism for $j=1$ by the same spectral sequence argument as in subsection \ref{Bdomain}.
By the remark after Definition \ref{defn:domain}, it is therefore a quasi-isomorphism for all $j$. Thus we 
have a quasi-isomorphism of domains $\delta : B_* \to W_*$.

Finally, it is immediate that $\delta$ is a map of partial $\mathcal{O}$-algebras since the $\mathcal{O}$-structures
are defined in the same way by the shuffle map and the left action of $\mathcal{O}$ on itself. It follows that $\delta$ is a quasi-isomorphism of partial $\mathcal{O}$-algebras.  It is routine to check, using the techniques already described, that all of our constructions were functorial, so this completes the proof of the theorem.

\section{Applications} \label{sec:appl}

In \cite{JM} McClure showed that the PL-chains $C$ of a PL-manifold have a
domain $\{C_j \}$ described by ``chains in general position''.
In particular, McClure showed that $C_j$ is quasi-isomorphic to $C^{\otimes_j}$.
Moreover, McClure showed this domain is part of a 
``partial Leinster algebra'', defined using the intersection of chains. In the language of this paper, 
this means they form a partial algebra over the operad 
$\mathscr{C}$ describing commutative associative algebras. 

A technical assumption of Theorem~\ref{thm:convert} is that the $j^{th}$  
component of the operad must be, for each $j > 0$, a projective $R[\Sigma_j]$-module.
Over $R=\mathbb{Q}$ every module is $R[\Sigma_j]$-projective, so 
by Theorem~\ref{thm:convert} we obtain from any partial $\mathscr{C}$-algebra $A$  (over $\mathbb{Q}$)  a commutative associative differential graded algebra on a complex quasi-isomorphic to $A$. In particular, we have the example $A=C$ above.

On the other hand, over $R = \mathbb{Z}$, $\mathscr{C}$ does not satisfy this 
property (since the $\Sigma_j$ actions are trivial on $\mathscr{C}(j) = 
\mathbb{Z}$). The following operad does satisfy the projective assumption:

\begin{defn}
An $\mathcal{E}_\infty$-operad is a unital operad $\mathcal{E}$, i.e. 
$\mathcal{E}(0) \approx \mathbb{Z}$, such that the maps
$$
\mathcal{E}(j) \otimes \mathcal{E}(0)^j \to \mathcal{E}(0) \approx \mathbb{Z}
$$ 
are quasi-isomorphisms, and each $\mathcal{E}(j)$ is a free 
$\mathbb{Z}[\Sigma_j]$-module. An algebra over an  $\mathcal{E}_\infty$-operad 
is called an  $\mathcal{E}_\infty$-algebra.
\end{defn}

Following Kriz and May in \cite{KM}, we use the given quasi-isomorphisms
$$
\mathcal{E}(j) \approx \mathbb{Z} = \mathscr{C}(j),
$$
to pull back the partial $\mathscr{C}$-algebra on the PL-chains $C$ of a PL-manifold, \cite{JM},
to obtain a partial algebra over $\mathcal{E}$. Then by 
Theorem~\ref{thm:convert} we obtain the following: 
 
\begin{thm}[McClure \cite{JM} using Theorem \ref{thm:convert} above] \label{thm:E}
There is a functor assigning to any closed PL-manifold an 
$\mathcal{E}_\infty$-algebra on a complex quasi-isomorphic to its PL-chains.
This structure induces the intersection product in homology.
\end{thm}

\begin{rmk}
It is natural to ask how this $\mathcal{E}_\infty$ chain-algebra relates to
the known $\mathcal{E}_\infty$ cochain-algebra that, by a theorem of Mandell
\cite{MM}, determines the weak homotopy type of a finite type nilpotent space\footnote{It has recently been announced by  D. Chataur that these two $\mathcal{E}_\infty$-algebras are in fact quasi-isomorphic.}.
\end{rmk}

We can take this example a bit further: the (PL) chains of a (PL) manifold 
embed quasi-isomorphically into the space of bounded
linear functionals on differential forms with compact support, i.e. currents, see de Rham \cite{DR}.
The same domain $\{C_j\}$ of chains as before also gives a domain for currents.
By the same argument as above, over $\mathbb{Q}$ or $\mathbb{R}$, Theorem~\ref{thm:convert} 
gives a commutative associative differential graded algebra on a 
complex quasi-isomorphic to the currents.
This gives an algebraic resolution to the long standing issue in functional 
analysis of not being able to multiply distributions, which are precisely zero currents.

\bigskip

{\sc Scott O. Wilson; 303 Kiely Hall; 65-30 Kissena, Blvd.; Flushing, NY 11367 USA.}

email: {\tt scott.wilson@qc.cuny.edu}

\end{document}